\documentclass[a4paper,reqno]{amsart}

\usepackage{amscd}
\usepackage{wrapfig}
\usepackage{graphicx}
\usepackage{epic,eepic}
\usepackage{hyperref}

\newcommand{\End}{\text{End}}
\newcommand{\uqgh}{U_q(\hat{g})}

\newcommand{\one}{\mathbf{1}}
\newcommand{\id}{\text{id}}
\newcommand{\Qb}{\overline{Q}}
\newcommand{\mub}{\bar{\mu}}
\newcommand{\nub}{\bar{\nu}}

\newcommand{\sae}{{\mathcal{B}_\epsilon}}
\newcommand{\ab}{\alpha}

\newcommand{\ur}{{\mathcal R}}

\newcommand{\aaa}{\uqgh}

\newcommand{\drawcenteredtext}[3]{\put(#1,#2){\makebox(0,0){#3}}}%
\newcommand{\drawlefttext}[3]{\put(#1,#2){\makebox(0,0)[l]{#3}}}%
\newcommand{\drawrighttext}[3]{\put(#1,#2){\makebox(0,0)[r]{#3}}}%

\newcommand{\drawpath}[4]{\path(#1,#2)(#3,#4)}%
\newcommand{\drawthickdot}[2]{\put(#1,#2){\circle*{1}}}%

\newcommand{\drawleftbrace}[3]%
{\drawcenteredtext{#1}{#2}{$\left\{ \rule[0mm]{0mm}{#3mm} \right.$}}%

\newcommand{\drawrightbrace}[3]%
{\drawcenteredtext{#1}{#2}{$\left\} \rule[0mm]{0mm}{#3mm} \right.$}}%

\newcommand{\drawoverbrace}[3]%
{\drawcenteredtext{#1}{#2}{$\overbrace{\rule[0mm]{#3mm}{0mm}}$}}%

\newcommand{\drawunderbrace}[3]%
{\drawcenteredtext{#1}{#2}{$\underbrace{\rule[0mm]{#3mm}{0mm}}$}}%

\newcommand{\drawarc}[5]%
{\put(#1,#2){\arc{#3}{#4}{#5}}}%

\begin{document}

\title[Solitons, Boundaries, and Quantum Affine Algebras]{Solitons, Boundaries,\\ and Quantum Affine Algebras}

\author{Gustav W Delius}

\address{Department of Mathematics\\University of
York\\York YO10 5DD\\United Kingdom}

\email{gwd2@york.ac.uk}

\urladdr{http://www.york.ac.uk/mathematics/physics/delius/}

\thanks{The
transparencies from this talk are available on the web at\\
\href{http://www.york.ac.uk/depts/maths/physics/delius/talks/iskmaa.pdf}{http://www.york.ac.uk/depts/maths/physics/delius/talks/iskmaa.pdf}
.}

\begin{abstract}

This is a condensed write-up of a talk delivered at the Ramanujan
International Symposium on Kac-Moody Lie algebras and Applications
in Chennai in January 2002. The talk introduces special coideal
subalgebras of quantum affine algebras which appear in physics
when solitons are restricted to live on a half-line by an
integrable boundary condition. We review how the quantum affine
symmetry determines the soliton S-matrix in affine Toda field
theory and then go on to use the unbroken coideal subalgebra on
the half-line to determine the soliton reflection matrix. This
gives a representation theoretic method for the solution of the
reflection equation (boundary Yang-Baxter equation) by reducing it
to a linear equation.
\end{abstract}

\maketitle

\section{Introduction}

{\bf Quantum affine algebras} are quantum deformations of the
universal enveloping algebras of affine Kac-Moody algebras,
introduced by Drinfeld \cite{Drinfeld:1988in} and Jimbo
\cite{Jimbo:1986ua}. They are therefore obviously related to the
topic of this conference.

{\bf Solitons} are localized finite-energy solutions of
relativistic wave equations with special properties. In the
quantum theory they lead to particle states. The relation between
classical solitons and Kac-Moody algebras was stressed already by
the Kyoto group \cite{MR2001a:37109}. More recently the role that
quantum affine algebras play in the quantum theory of solitons was
uncovered in \cite{Bernard:1991ys}.

By {\bf Boundaries} I mean spatial boundaries. In this talk space
will be taken to be one-dimensional, i.e., the real line.
Introducing a boundary by imposing a boundary condition at some
point restricts physics to the half-line.

In this talk we will knit these three topics together. We will do
so by using a particularly rich family of quantum field theories.
They are known as quantum affine Toda field theories and they
posses all three: a quantum affine symmetry algebra
\cite{Bernard:1991ys}, soliton solutions
\cite{Hollowood:1992by,Olive:1993cm}, and integrable boundary
conditions \cite{Bowcock:1995vp}.

The main topic of the talk will be certain coideal subalgebras of
quantum affine algebras that remain unbroken symmetry algebras
after imposing an integrable boundary condition
\cite{Delius:2001qh}. These algebras provide a new tool for the
solution of the reflection equation.

Because the first part of the talk is simply an introduction to
the known theory of solitons and quantum affine symmetry in affine
Toda field theory, we will be very brief in the next two sections.
We give references to the literature where more details can be
found. These references are not always to the original works but
often to more recent works because they may be more useful in the
given context.

\section{Solitons in Affine Toda field theory}\label{s:toda}

\subsection{Field equations}
Associated to every affine Kac-Moody algebra $\hat{g}$ there is a relativistic field
equation
\begin{equation}\label{feq}
  \partial_x^2\phi-\partial_t^2\phi=\frac{1}{2i}\sum_{j=0}^n\,\eta_j\,\alpha_j^\vee\,
  e^{i(\alpha_j^\vee,\,\phi)}
\end{equation}
where the field $\phi=\phi(x,t)$ takes values in the root space of
the underlying finite dimensional Lie algebra $g$,
$\alpha_0,\dots,\alpha_n$ are the simple roots of  $\hat{g}$
projected onto the root space of $g$, $(\ ,\ )$ is the Killing
form, $\alpha_j^\vee=2\alpha_j/(\alpha_j,\alpha_j)$ and the
$\eta_j$ are coprime integers such that
$\sum_{j=0}^n\eta_j\alpha_j^\vee=0$. The sine-Gordon model is the
simplest example of an affine Toda field equation corresponding to
$g=\hat{sl}_2$.

The field equations \eqref{feq} are integrable in the sense that
they posses an infinite number of integrals of motion in
involution. This follows for example from the fact that they can
be written as the zero-curvature condition for a gauge-connection
taking its value in the affine Lie algebra $\hat{g}$
\cite{Olive:1985mb}.

\subsection{Soliton solutions}
The Toda field equations \eqref{feq} have degenerate vacuum
solutions with $\phi=2\pi\lambda$ for any fundamental weight
$\lambda$ of $g$. Thus there are kink solutions which interpolate
between these vacua. They are stable due to their nonzero
topological charge
\begin{equation}\label{tcharge}
  T[\phi]=\phi(\infty)-\phi(-\infty)=2\pi\lambda.
\end{equation}

\begin{wrapfigure}{r}{4cm}
\includegraphics[width=4cm]{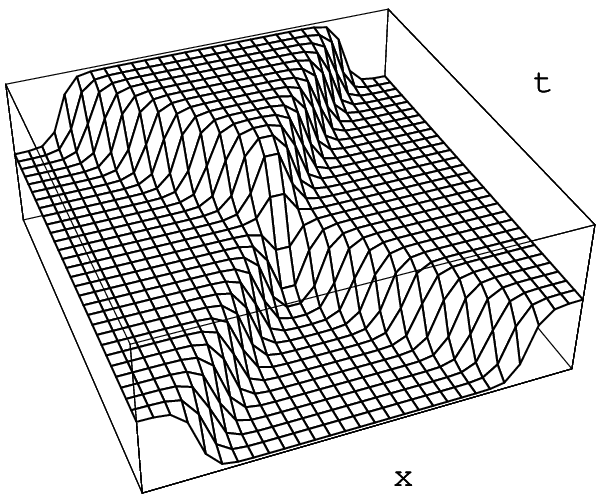}
\end{wrapfigure}
Furthermore there are solutions that describe an arbitrary number
of kinks, each moving with its own velocity. Surprisingly, when
two of these kinks meet they pass through each other and reemerge
with their original shape restored. It is  this property that
shows that these kinks are solitons. An example of a two-soliton
solution in the sine-Gordon model is shown on the right.

The topological charges of the fundamental solitons are weights in
the fundamental representations of $g$. All solitons whose
topological charge lie in the same fundamental representation have
the same mass. A particle physicist would therefore say that they
form a multiplet.

As an introduction to solitons from the particle physicsts
viewpoint I recommend the book by Rajaraman
\cite{Rajaraman:1982is}. A more classical presentation is given in
\cite{MR2001a:37109}. For details on the solitons in affine Toda
theory see \cite{Hollowood:1992by,Olive:1993cm}.

\subsection{Quantum solitons}

In the quantum theory we associate particle states with the
soliton solutions. Let $V^\mu_\theta$ be the space spanned by the
solitons in multiplet $\mu$ with rapidity\footnote{The rapidity
$\theta$ parametrizes the energy and momentum of particles on the
mass shell. If the mass is $m$ then the energy is $E=m\cosh\theta$
and the momentum is $p=m\sinh\theta$.} $\theta$.

\begin{wrapfigure}{r}{4cm}
\setlength{\unitlength}{1.0mm}
\begin{picture}(66,15)(15,9)
\thinlines \drawpath{22.0}{22.0}{34.0}{10.0}
\drawpath{34.0}{22.0}{22.0}{10.0} \drawthickdot{28.0}{16.0}
\drawcenteredtext{22.0}{6.0}{$V^\mu_\theta$}
\drawcenteredtext{34.0}{6.0}{$V^\nu_{\theta'}$}
\drawcenteredtext{20.0}{26.0}{$V^\nu_{\theta'}$}
\drawcenteredtext{34.0}{26.0}{$V^\mu_\theta$}
\drawlefttext{34.0}{16.0}{$S^{\mu\nu}(\theta-\theta')$}
\end{picture}
\end{wrapfigure}

Asymptotic two-soliton states span tensor product spaces
$V^\mu_\theta\otimes V^\nu_{\theta'}$.

An incoming two-soliton state in $V^\mu_\theta\otimes
V^\nu_{\theta'}$ with $\theta>\theta'$ will evolve during
scattering into an outgoing state in $V^\nu_{\theta'}\otimes
V^\mu_{\theta}$ with scattering amplitude given by the entry of
the two-soliton S-matrix $S^{\mu\nu}(\theta-\theta')$. We
represent this graphically as in the figure on the right.

\subsection{S-matrix factorization}

Due to integrability, the multi-soliton S-matrix factorizes into a
product of two-soliton S-matrices. For the three-particle
scattering process this is represented graphically as follows:
\begin{equation}
\setlength{\unitlength}{0.8mm}
\begin{picture}(112,32)
\thinlines \drawpath{4.0}{28.0}{28.0}{4.0}
\drawpath{28.0}{28.0}{4.0}{4.0} \drawpath{16.0}{28.0}{16.0}{4.0}
\drawpath{44.0}{28.0}{68.0}{4.0} \drawpath{68.0}{28.0}{44.0}{4.0}
\drawpath{62.0}{28.0}{62.0}{4.0} \drawpath{84.0}{28.0}{108.0}{4.0}
\drawpath{108.0}{28.0}{84.0}{4.0} \drawpath{92.0}{28.0}{92.0}{4.0}
\drawcenteredtext{38.0}{16.0}{$=$}
\drawcenteredtext{78.0}{16.0}{$=$} \drawthickdot{16.0}{16.0}
\drawthickdot{56.0}{16.0} \drawthickdot{62.0}{22.0}
\drawthickdot{62.0}{10.0} \drawthickdot{92.0}{20.0}
\drawthickdot{96.0}{16.0} \drawthickdot{92.0}{12.0}
\end{picture}
\end{equation}
The compatibility of the two ways to factorize requires the
S-matrix to be a solution of the Yang-Baxter equation
\begin{multline}\label{ybe}
    (1\otimes S^{\mu\nu}(\theta-\theta'))
    (S^{\mu\lambda}(\theta-\theta'')\otimes 1)
    (1\otimes S^{\nu\lambda}(\theta'-\theta''))
    \\=
    (S^{\nu\lambda}(\theta'-\theta'')\otimes 1)
    (1\otimes S^{\mu\lambda}(\theta-\theta''))
    (S^{\mu\nu}(\theta-\theta')\otimes 1).
\end{multline}
A classical reference on factorized S-matrices is
\cite{Zamolodchikov:1979xm}.

Quantum affine algebras provide the key technique for finding
solutions to the Yang-Baxter equation, as will be explained in the
next section.

\section{Quantum group symmetry}

\subsection{Non-local symmetry charges}
Quantum affine Toda theory has symmetry charges $T_i, Q_i,
\bar{Q}_i$, $i=0,\dots,n$ which generate the $U_q(\hat{g})$
algebra with relations
\begin{align}
  &[T_i,Q_j]=\ab_i\cdot\ab_j\,Q_j,~~~~~~
  [T_i,\bar{Q}_j]=-\ab_i\cdot\ab_j\,Q_j\notag\\
  &Q_i\bar{Q}_j-q^{-\alpha_i\cdot\alpha_j}\bar{Q}_j
  Q_i=\delta_{ij}\,\frac{q^{2 T_i}-1}{q_i^{2}-1},
\end{align}
where
\begin{equation}
  q=e^{2\pi i\frac{1-\hbar}{\hbar}}
\end{equation}
and $q_i=q^{\ab_i\cdot\ab_i/2}$. They also satisfy the Serre
relations. The charges $Q_i, \bar{Q}_i$ are non-local in the sense
that they are obtained as space integrals of the time components
of currents which are themselves non-local expressions. The
details can be found in \cite{Bernard:1991ys}.

\subsection{Action on solitons}
Each soliton multiplet $V^\mu_\theta$ carries a representation
$\pi^\mu_\theta: {\aaa}\rightarrow \End(V^\mu_\theta)$ of the
quantum affine algebra. The symmetry acts on the multi-soliton
states through the coproduct $\Delta:{\aaa}\rightarrow
{\aaa}\otimes {\aaa}$.
\begin{align}\label{cop}
  \Delta(Q_i)&=Q_i\otimes 1+q^{T_i}\otimes Q_i,\notag\\
  \Delta(\Qb_i)&=\Qb_i\otimes 1+q^{T_i}\otimes \Qb_i,\\
  \Delta(T_i)&=T_i\otimes 1+1\otimes T_i.\notag
\end{align}
An explanation of why non-local symmetry charges have such a
non-cocommutative coproduct can be found in \cite{Bernard:1993mu}.

\subsection{S-matrix as intertwiner}
The S-matrix has to commute with the action of any symmetry charge
$Q\in{\aaa}$,
\begin{equation}
\begin{CD}\label{sint}
  V^\mu_{\theta}\otimes V^\nu_{\theta'}
  @>(\pi^\mu_\theta\otimes\pi^\nu_{\theta'})(\Delta(Q))>>
  V^\mu_{\theta}\otimes V^\nu_{\theta'}
  \\
  @VV{S^{\mu\nu}(\theta-\theta')}V
  @VV{S^{\mu\nu}(\theta-\theta')}V
  \\
  V^\nu_{\theta'}\otimes V^\mu_{\theta}
  @>(\pi^\nu_{\theta'}\otimes\pi^\mu_{\theta})(\Delta(Q))>>
  V^\nu_{\theta'}\otimes V^\mu_{\theta}
\end{CD}
\end{equation}
This determines the S-matrix uniquely up to an overall factor.
Jimbo solved this intertwining equation for the vector
representations of all quantum affine algebras in
\cite{Jimbo:1986ua}. A practical technique for solving the
intertwining equations in a large number of cases is the tensor
product graph method, see e.g. \cite{Delius:1994ht}.

The scalar prefactor which is not fixed by the Yang-Baxter
equation can be determined by imposing other physical requirements
such as unitarity, crossing symmetry and closure of the bootstrap.
This has been done for a number of affine Toda theories, see e.g.,
\cite{Hollowood:1993sy,Gandenberger:1995gg,Gandenberger:1996cw}

\subsection{Yang-Baxter equation from Schur's lemma}
Because the tensor product modules are irreducible for generic
rapidities, Schur's lemma implies that the following diagram is
commutative up to a scalar factor\begin{equation}
\begin{CD}
  V^\mu_{\theta}\otimes V^\nu_{\theta'}\otimes V^\lambda_{\theta''}
  @>{S^{\mu\nu}(\theta-\theta')\,\otimes\,\id}>>
  V^\nu_{\theta'}\otimes V^\mu_{\theta}\otimes V^\lambda_{\theta''}
  \\
  @VV{\id\otimes S^{\nu\lambda}(\theta'-\theta'')}V
  @V{id\,\otimes\, S^{\mu\lambda}(\theta-\theta'')}VV
  \\
  V^\mu_{\theta}\otimes V^\lambda_{\theta''}\otimes V^\nu_{\theta'}
  @.
  V^\nu_{\theta'}\otimes V^\lambda_{\theta''}\otimes V^\mu_{\theta}
  \\
  @VV{S^{\mu\lambda}(\theta-\theta'')\,\otimes\,\id}V
  @V{S^{\nu\lambda}(\theta'-\theta'')\otimes\id}VV
  \\
  V^\lambda_{\theta''}\otimes V^\mu_{\theta}\otimes V^\nu_{\theta'}
  @>{id\,\otimes\, S^{\mu\nu}(\theta-\theta')}>>
  V^\lambda_{\theta''}\otimes V^\nu_{\theta'}\otimes V^\mu_{\theta}
\end{CD}
\end{equation}
Thus any S-matrix which satisfies the intertwining property
\eqref{sint} automatically satisfies the Yang-Baxter equation (at
least up to a scalar factor).

\section{Soliton reflection}

\subsection{Integrable boundary conditions}
We now restrict the theory to the left half-line by imposing that
at the boundary
\begin{equation}\label{bc}
\partial_x{\phi} = i\sum_{j=0}^n \epsilon_j
\alpha_j^\vee\exp\left(\frac{i}{2}(\alpha_j^\vee,\phi)\right)
\end{equation}
where the $\epsilon_j$ are free parameters. These boundary
conditions were proposed by Corrigan et.al. because they preserve
integrability for certain choices of the parameters
\cite{Bowcock:1995vp}. There are other integrable boundary
conditions \cite{Bowcock:1996gw,Delius:1998rf} and it will be
interesting to extend the observations from the next section to
these.

\subsection{Soliton reflection}

\begin{wrapfigure}{r}{3.5cm}
\setlength{\unitlength}{0.8mm}
\begin{picture}(62,21)(18,10)
\Thicklines \drawpath{40.0}{26.0}{40.0}{6.0} \thinlines
\drawpath{40.0}{16.0}{30.0}{26.0} \drawpath{40.0}{16.0}{30.0}{6.0}
\drawrighttext{28.0}{6.0}{$V^\mu_\theta$}
\drawrighttext{28.0}{26.0}{$V^{\bar{\mu}}_{-\theta}$}
\drawlefttext{42.0}{16.0}{$K^\mu(\theta)$}
\end{picture}
\end{wrapfigure}
For $a_n^{(1)}$ affine Toda theory the soliton solutions which
satisfy the boundary conditions \eqref{bc} were determined in
\cite{Delius:1998jw} by a nonlinear method of images. These
solutions describe an incoming soliton with rapidity $\theta$
being converted into an outgoing anti-soliton with rapidity
$-\theta$. In the quantum theory such reflection processes are
described by a reflection matrix $K^\mu(\theta)$ which maps from
the incoming multiplet $V^\mu_\theta$ to the outgoing conjugate
multiplet $V^{\bar{\mu}}_{-\theta}$.
The entries of the reflection
matrix are the reflection amplitudes.

\subsection{Factorization of reflection matrices}
If the boundary condition preserves integrability then the
multi-soliton reflection matrices have to factorize.
\begin{equation}
\setlength{\unitlength}{0.8mm}
\begin{picture}(98,38)
\Thicklines \drawpath{20.0}{34.0}{20.0}{4.0} \thinlines
\drawpath{20.0}{20.0}{6.0}{34.0} \drawpath{20.0}{20.0}{6.0}{6.0}
\drawpath{20.0}{20.0}{4.0}{28.0} \drawpath{20.0}{20.0}{4.0}{12.0}
\drawthickdot{20.0}{20.0} \drawcenteredtext{30.0}{20.0}{$=$}
\drawcenteredtext{66.0}{20.0}{$=$} \Thicklines
\drawpath{56.0}{34.0}{56.0}{4.0} \thinlines
\drawpath{56.0}{24.0}{36.0}{4.0} \drawpath{56.0}{24.0}{46.0}{34.0}
\drawpath{56.0}{18.0}{36.0}{28.0} \drawpath{56.0}{18.0}{36.0}{8.0}
\drawthickdot{52.0}{20.0} \drawthickdot{44.0}{12.0}
\drawthickdot{56.0}{18.0} \drawthickdot{56.0}{24.0} \Thicklines
\drawpath{92.0}{34.0}{92.0}{4.0} \thinlines
\drawpath{92.0}{14.0}{82.0}{4.0} \drawpath{92.0}{14.0}{72.0}{34.0}
\drawpath{92.0}{20.0}{72.0}{30.0}
\drawpath{92.0}{20.0}{72.0}{10.0} \drawthickdot{80.0}{26.0}
\drawthickdot{92.0}{20.0} \drawthickdot{88.0}{18.0}
\drawthickdot{92.0}{14.0}
\end{picture}
\end{equation}
The compatibility condition between the two ways of factorizing
the two-soliton reflection matrix is the {\it reflection equation}
\cite{Cherednik:1984vs}
\begin{multline}\label{req}
    (1\otimes K^\nu(\theta'))
    S^{\nu\bar{\mu}}(\theta+\theta')
    (1\otimes K^\mu(\theta))
    S^{\mu\nu}(\theta-\theta')
    \\=
    S^{\bar{\nu}\bar{\mu}}(\theta-\theta')
    (1\otimes K^\mu(\theta))
    S^{\mu\bar{\nu}}(\theta+\theta')
    (1\otimes K^\nu(\theta').
\end{multline}
This equation is also sometimes known as the {\it boundary
Yang-Baxter equation}. Like the Yang-Baxter equation \eqref{ybe}
it is a non-linear functional matrix equation. It was solved to
find the reflection matrix for the sine-Gordon solitons in
\cite{Ghoshal:1994tm} and then for the vector solitons in
$a_n^{(1)}$ affine Toda theory in \cite{Gandenberger:1999uw}.

Solutions of the reflection equation are not only needed to
describe the reflection of particles off a boundary. They are also
used to construct integrable quantum spin chains with a boundary
and they provide the boundary Bolzmann weights for solvable models
of statistical mechanics. Unfortunately the reflection equation
has proven difficult to solve directly except for matrices of very
small dimension \cite{Lima-Santos:1999,Liu:1998jd}, for diagonal
matrices \cite{deVega:1994sb,Batchelor:1996xa} and a few other
cases \cite{Ahn:1998,Gandenberger:1999uw,Lima-Santos:2001aa}. The
two known systematic techniques for the construction of reflection
matrices are the fusion procedure \cite{Behrend:1996en} and the
extended Hecke algebraic techniques \cite{Doikou:2002ry} but these
have so far only produced a small subset of the possible
solutions.

A way to find solutions to the reflection equation with the help
of quantum affine algebras has only been discovered very recently
\cite{Delius:2001qh}. This method arose by studying the non-local
symmetries of affine Toda theory on the half-line, as will be
described in the next section.

\section{Boundary quantum group symmetry}

\subsection{Symmetry charges}
In the presence of the boundary the non-local charges of the bulk
theory are no longer conserved. Rather we found
\cite{Delius:2001qh} that the boundary condition \eqref{bc} breaks
the symmetry to a subalgebra ${\mathcal B}_\epsilon\subset{\aaa}$
generated by
\begin{equation}\label{tcc}
  \widehat{Q}_i = Q_i + \bar{Q}_i +
  \hat{\epsilon}_i q^{T_i},~~~~~i=0,\dots,n.
\end{equation}
Note that the symmetry algebra depends on the boundary parameters
$\hat{\epsilon_i}$ (which are related to the $\epsilon_i$ in
\eqref{bc} by a rescaling). We derived these conserved charges
using first order boundary conformal perturbation theory but we
believe the result to hold exactly. For $\hat{g}=a_n^{(1)}$ we
were able to show that this is the case, at least on shell, as
explained in the next section. We can now report that recent
calculations performed in collaboration with Alan George establish
the symmetry also for $\hat{g}=d_n^{(1)}$. It is the subject of
ongoing work to extend this to all affine Toda field theories with
an integrable boundary condition.

\subsection{Reflection Matrix as Intertwiner}

The reflection matrix has to commute with the action of any
symmetry charge $\hat{Q}\in{\mathcal B}_\epsilon\subset{\aaa}$,
i.e., the following diagram is commutative:
\begin{equation}\label{kint}
\begin{CD}
  V^\mu_{\theta}
  @>\pi^\mu_\theta(\hat{Q})>>
  V^\mu_{\theta}
  \\
  @VV{K^{\mu}(\theta)}V
  @VV{K^{\mu}(\theta)}V
  \\
  V^{\mub}_{-\theta}
  @>\pi^{\mub}_{-\theta}(\hat{Q})>>
  V^{\mub}_{-\theta}
\end{CD}
\end{equation}
If the symmetry algebra $\sae$ is "large enough" so that
$V^\mu_\theta$ and $V^{\mub}_{-\theta}$ are irreducible modules of
$\sae$, then this linear equation determines the reflection matrix
uniquely up to an overall factor. As we will see below, this is
the case for example for the vector representation of $a_n^{(1)}$.

\subsection{Coideal property}

The residual symmetry algebra $\sae$ turns out not to be a Hopf
algebra because it does not have a coproduct
$\Delta:\sae\rightarrow\sae\otimes\sae$. Using the coproduct
\eqref{cop} of $\aaa$ we find
\begin{equation}
  \Delta(\hat{Q}_i)=(Q_i+\bar{Q}_i)\otimes 1+q^{T_i}\otimes
  \hat{Q}_i.
\end{equation}
Thus the symmetry algebra $\sae$ is a left coideal of ${\aaa}$ in
the sense that
\begin{equation}
  \Delta(\hat{Q})\in{\aaa}\otimes\sae~~~\text{ for all }\hat{Q}\in\sae.
\end{equation}
This allows it to act on multi-soliton states in the presence of a
boundary because the coproduct can be used to define the action of
$\sae$ on any tensor product of a $\aaa$-module (solitons) with a
$\sae$-module (boundary states). We would like to stress that this
is the first time that a symmetry algebra in physics has been
found that is not a bialgebra.

\subsection{The Reflection Equation from Schur's lemma}
If ${\mathcal B}_\epsilon$ is "large enough" so that the tensor
product modules are irreducible, then by Schur's lemma the
following diagram is commutative up to a scalar factor
\begin{equation}
\begin{CD}\label{cdre}
  V^\mu_\theta\otimes V^\nu_{\theta'}
  @>\id\,\otimes K^\nu(\theta')>>
  V^\mu_\theta\otimes V^{\bar{\nu}}_{-\theta'}
  \\
  @VVS^{\mu\nu}(\theta-\theta')V
  @VS^{\mu\bar{\nu}}(\theta+\theta')VV
  \\
  V^\nu_{\theta'}\otimes V^\mu_{\theta}
  @.
  V^{\bar{\nu}}_{-\theta'}\otimes V^\mu_\theta
  \\
  @VV\id\,\otimes K^\mu(\theta)V
  @VV\id\,\otimes K^\mu(\theta)V
  \\
  V^\nu_{\theta'}\otimes V^{\bar{\mu}}_{-\theta}
  @.
  V^{\bar{\nu}}_{-\theta'}\otimes V^{\bar{\mu}}_{-\theta}
  \\
  @VVS^{\nu\bar{\mu}}(\theta+\theta')V
  @VS^{\bar{\nu}\bar{\mu}}(\theta-\theta')VV
  \\
  V^{\bar{\mu}}_{-\theta}\otimes V^\nu_{\theta'}
  @>\id\,\otimes K^\nu(\theta')>>
  V^{\bar{\mu}}_{-\theta}\otimes V^{\bar{\nu}}_{-\theta'}
\end{CD}
\end{equation}
Thus any reflection matrix which satisfies the intertwining
property \eqref{kint} automatically satisfies the reflection
equation \eqref{req} (at least up to a scalar factor).

\subsection{Calculating Reflection Matrices}
Solving the quantum group intertwining property \eqref{kint} is
much easier than solving the reflection equation \eqref{req}
directly because the intertwining property is a linear equation.

To illustrate this new technique we present here the derivation of
the reflection matrix for the vector solitons in $a_n^{(1)}$ Toda
theory. Using the representation matrices
\begin{align}
  \pi^\mu_\theta(\hat{Q}_i)&=x\,e^{i+1}{}_i+
  x^{-1}\,e^{i}{}_{i+1}+
  \hat{\epsilon}_i\,((q^{-1}-1)\,e^{i}{}_i+(q-1)\,e^{i+1}{}_{i+1}+\one)
\end{align}
the intertwining property \eqref{kint} gives the following set of
linear equations for the entries of the reflection matrix:
\begin{align}
  0&=\hat{\epsilon}_i(q^{-1}-q)K^i{}_i+x\,K^i{}_{i+1}-x^{-1}\,K^{i+1}{}_i,
  \label{s3}\\
  0&=K^{i+1}{}_{i+1}-K^i{}_i,\\
  0&=\hat{\epsilon}_i\,q\,K^i{}_j+x^{-1}\,K^{i+1}{}_j,~~~j\neq
  i,i+1,\\
  0&=\hat{\epsilon}_i\,q^{-1}\,K^j{}_i+x\,K^j{}_{i+1},~~~j\neq
  i,i+1.
\end{align}

If all $|\hat{\epsilon}_i|=1$ then one finds the solution
\begin{align}
  K^i{}_i(\theta)&=\left(
  q^{-1}\,(-q\,x)^{(n+1)/2}-
  \hat{\epsilon}\,q\,(-q\,x)^{-(n+1)/2}
  \right)
  \frac{k(\theta)}{q^{-1}-q},\\
  K^i{}_j(\theta)&=\hat{\epsilon}_i\cdots\hat{\epsilon}_{j-1}\,
  (-q\,x)^{i-j+(n+1)/2}\,k(\theta),~~~~~~\text{for }j>i,\\
  K^j{}_i(\theta)&=\hat{\epsilon}_i\cdots\hat{\epsilon}_{j-1}\hat{\epsilon}\,
  (-q\,x)^{j-i-(n+1)/2}\,k(\theta),~~~~~\text{for }j>i,
\end{align}
that is unique up to an overall numerical factor $k(\theta)$. If
all $\hat{\epsilon}_i=0$ then the solution is diagonal. For other
values for the $\hat{\epsilon}_i$ there are no solutions. This is
in agreement with the restrictions on the boundary parameters
found by imposing classical integrability \cite{Bowcock:1995vp}.

We have applied the same technique to derive the reflection
matrices for any representation of $\hat{sl}_2$
\cite{Delius:2002ad}. We are currently calculating the reflection
matrices for the vector representation for all families of affine
Lie algebras, mirroring the work done by Jimbo for R-matrices
\cite{Jimbo:1986ua}.

\subsection{Boundary Bound States}
Particles can bind to the boundary, creating multiplets of
boundary bound states. These span finite dimensional
representations $V^{[\lambda]}$ of the symmetry algebra ${\mathcal
B}_\epsilon$. The reflection of particles off boundary bound
states is described by intertwiners
\begin{equation}
  K^{\mu[\lambda]}(\theta):V^\mu_\theta\otimes
V^{[\lambda]}\rightarrow V^{\mub}_{-\theta}\otimes V^{[\lambda]}.
\end{equation}
Unfortunately nothing is known so far about the representation
theory of the algebras $\sae$. One objective of this talk is to
encourage research into this direction. Ideally one would like to
be able to calculate the branching rules for representations of
$\aaa$ into irreducible representations of $\sae$.

\section{Reconstruction of symmetry}
In this last section I want to show how to reconstruct the
symmetry algebra of a theory with boundary from the knowledge of
at least one of its reflection matrices.

Let us assume that for one particular $\aaa$ representation
$V^\mu_\theta$ we know the reflection matrix
$K^\mu(\theta):V^\mu_\theta\rightarrow V^{\bar{\mu}}_{-\theta}$.
We define the corresponding ${\aaa}$-valued $L$-operators in terms
of the universal $R$-matrix ${\mathcal R}$ of ${\aaa}$,
\begin{align}
  L^\mu_\theta&=\left(\pi^\mu_\theta\otimes\id\right)
  ({\mathcal R})\in\End(V^\mu_\theta)\otimes{\aaa},\\
  \bar{L}^{\bar{\mu}}_\theta&=
  \left(\pi^{\bar{\mu}}_{-\theta}\otimes\id\right)
  ({\mathcal R}^{\text{op}})\in\End(V^{\mub}_{-\theta})\otimes{\aaa}.
\end{align}
Here $\ur^{\text{op}}$ is the opposite universal
R-matrix obtained by interchanging the two tensor factors.
From these $L$-operators we construct the matrices
\begin{equation}
  B^\mu_\theta=\bar{L}^{\bar{\mu}}_\theta\,
  (K^\mu(\theta)\otimes\,1)\,L^\mu_\theta\in
  \End(V^\mu_\theta,V^{\mub}_{-\theta})\otimes{\aaa}.
\end{equation}
In index notation this becomes
\begin{equation}
  (B^\mu_\theta)^\alpha{}_\beta=
  (\bar{L}^{\bar{\mu}}_\theta)^\alpha{}_\gamma
  (K^\mu(\theta))^\gamma{}_\delta(L^\mu_\theta)^\delta{}_\beta
  \in{\aaa}.
\end{equation}
We find that for all $\theta$ the $(B^\mu_\theta)^\alpha{}_\beta$
are elements of a coideal subalgebra ${\mathcal B}$ which commutes
with the reflection matrices. It is easy to check the coideal
property:
\begin{equation}
  \Delta\left((B^{\mu}_\theta)^\alpha{}_\beta\right)
  =(\bar{L}^{\mub}_\theta)^\alpha{}_\delta
  (L^\mu_\theta)^\sigma{}_\beta
  \otimes (B^{\mu}_\theta)^\delta{}_\sigma.
\end{equation}
Also any $K^\nu(\theta'):V^\nu_{\theta'}\rightarrow
V^{\nub}_{-\theta'}$ which satisfies the appropriate reflection
equation commutes with the action of the elements
$(B^\mu_\theta)^\alpha{}_\beta$
\begin{equation}
  K^\nu(\theta')\circ\pi^\nu_{\theta'}((B^\mu_\theta)^\alpha{}_\beta)=
  \pi^{\bar{\nu}}_{-\theta'}((B^\mu_\theta)^\alpha{}_\beta)\circ K^\nu(\theta'),
\end{equation}

Applying the above construction to the vector solitons in affine
Toda theory and expanding in powers of $x=e^\theta$ gives
\begin{equation}
  B^\mu_\theta=B+x\,\sum_{l=0}^n(q^{-1}-q)\,e^{l+1}{}_l\otimes
  \left(Q_l+\bar{Q}_l+\hat{\epsilon}_l\,q^{T_l}\right)+{\mathcal O}(x^2).
\end{equation}
This shows that the charges were correct to all orders.

The $B$-matrices can be shown to satisfy the quadratic relations
\begin{equation}
  \check{P}R^{\nub\mub}(\theta-\theta')\check{P}\,
  \overset{1}{B^\mu_\theta}\,
  R^{\mu\nub}(\theta+\theta')\,
  \overset{2}{B^\nu_{\theta'}}
  =\overset{2}{B^\nu_{\theta'}}\,
  \check{P}R^{\nu\mub}(\theta+\theta')\check{P}\,
  \overset{1}{B^\mu_\theta}\,
  R^{\mu\nu}(\theta-\theta'),
\end{equation}
Thus they generate a {\it reflection equation algebra} in the
sense of Sklyanin \cite{Sklyanin:1988yz}.

\section{Summary and Discussion}

In the first half of the talk we reviewed how quantum affine
algebras are used to find the solutions of the Yang-Baxter
equation that describe soliton scattering in affine Toda theory.
In the second half we then described the recently discovered
analogous technique for soliton reflection off integrable
boundaries. The main points were the following
\begin{itemize}

\item A boundary condition breaks the symmetry to a subalgebra
$\sae$ of the quantum affie algebra $\aaa$, depending on the
boundary parameters $\epsilon_i$.

\item The residal symmetry algebra $\sae$ is not a Hopf algebra
but only a coideal subalgebra of $\aaa$.

\item The symmetry implies that the reflection matrices are
intertwiners of $\sae$ representations. Because for generic
rapidity the $\aaa$ modules $V^\mu_\theta$ spanned by the solitons
are irreducible also as modules of the subalgebra $\sae$, the
intertwining property determines the reflection matrices uniquely
up to a scalar factor.

\item Because also the tensor product modules $V^\mu_\theta\otimes
V^\nu_{\theta'}$ are irreducible as $\sae$ modules, Schur's lemma
implies that the reflection matrices that satisfy the intertwining
property are also automatically solutions of the reflection
equation.

\item The intertwining property leads to a set of linear equations
for the entries of the reflection matrix which is relatively
straightforward to solve. We have thus obtained a very practical
method of finding solutions to the reflection equation.

\item At special imaginary values of the rapidity $\theta$ a
soliton multiplets $V^\mu_theta$ can become reducible and boundary
bound states can form, spanning irreducible representations of
$\sae$. Thus the physically interesting study of boundary bound
states is connected to the study of the branching rules from
$\aaa$ to $\sae$.

\end{itemize}

While in this talk we concentrated on affine Toda field theories,
a similar story unfolds in the principal chiral models. There the
theory in the bulk has a Yangian symmetry, which in the presence
of a boundary is broken to a twisted Yangian. This was discovered
in \cite{Delius:2001yi} and was used to find reflection matrices
in \cite{Delius:2001yi,MacKay:2002at}.

We found the boundary quantum groups $\sae$ by studying concrete
boundary quantum field theories. However the question can also be
approached purely mathematially. One would like to classify and
construct all coideal subalgebras of quantum affine algebras or
Yangians that qualify as boundary quantum groups in the sense that
their intertwiners are solutions to the boundary Yang-Baxter
equation. This would lead to a classification of solutions of the
boundary Yang-Baxter equation in terms of boundary quantum groups
in the same spirit as the classification of solutions of the
Yang-Baxter equation in terms of representations of quantum
groups.

I hope to have given you in this talk a glimpse of a whole new
field of enquiry related to Kac-Moody algebras with applications
in both physics and mathematics.

\noindent{\bf Note added:}

Since this contribution was written in May 2002 several
interesting papers have appeared on the subject of reflection
matrices and coideal subalgebras. Examples include
\cite{Baseilhac:2002kf,Molev:2002,Delius:2002}

\noindent{\bf Acknowledgments}

The work presented here was performed in collaboration with Niall
MacKay and our graduate students Ben Short and Alan George. My
work is supported by an EPSRC advanced fellowship.

\providecommand{\href}[2]{#2}\begingroup\raggedright\endgroup

\end{document}